\documentclass{amsart}
\usepackage{amsmath, amsthm, amscd, amsfonts, amssymb, graphicx, color}
\usepackage[bookmarksnumbered]{hyperref}
\newtheorem{thm}{Theorem}[section]
\newtheorem{cor}[thm]{Corollary}

\theoremstyle{definition}
\newtheorem{defn}[thm]{Definition}
\theoremstyle{remark}

\numberwithin{equation}{section}

\begin{document}

\title[On a type of exponential
functional equation]{On a type of exponential
functional equation and its superstability in the sense of Ger}%
\author{A. Sousaraei $^{1}$, M. Alimohammady $^{2}$ and A. Sadeghi $^{3}$}%
\address{$^{1}$ Islamic Azad University Branch of Azadshaher, Azadshaher,
Iran.}
\address{$^{2}$ Department of Mathematics, University of Mazandaran,
Babolsar, Iran.}
\address{$^{3}$ Department of Mathematics Tarbiat Modares University Tehran 14115-333 Iran}%
\email{$^{1}$ ali Sousaraei@yahoo.com, $^{2}$ amohsen@umz.ac.ir, $^{3}$ ali.sadeghi@modares.ac.ir}%
\thanks{ }%
\subjclass[2000]{Primary 39B72, 39B82.}
\keywords{ Exponential functional equation, stability, superstability }%

\begin{abstract}
In this paper, we deal with a type exponential functional equation
as follows
$$f(xy)=f(x)^{g(y)},$$
where $f$ and $g$ are two real valued functions on a commutative
semigroup. Our aim of this paper is to proved that the above
functional equation in the sense of Ger is superstable.
\end{abstract}
\maketitle
\section{Introduction}

In 1979, a type of stability was observed by J. Baker, J. Lawrence
and F. Zorzitto \cite{b3}. Indeed, they proved that if a function
is approximately exponential, then it is either a true exponential
function or bounded. Then the exponential functional equation is
said to be superstable. This result was the first result
concerning the superstability phenomenon of functional equations.
Later, J. Baker \cite{b2} (see also \cite{17,Ali,asd,zxc})
generalized this famous result as follows:

Let $(S,\cdot)$ be an arbitrary semigroup, and let $f$ map $S$ into
the field $C$ of all complex numbers. Assume that $f$ is an
approximately exponential function, i.e., there exists a nonnegative
number $\varepsilon$ such that $$\|f(x\cdot
y)-f(x)f(y)\|\leq\varepsilon$$ for all  $x, y\in S$. Then $f$ is
either bounded or exponential.

The result of Baker, Lawrence and Zorzitto \cite{b3} was generalized
by L. Sz$\acute{e}$kelyhidi \cite{sz1} in another way and he
obtained the following result.
\begin{thm}\label{theo1}
Let $(G,.)$ be an Abelian group with identity $1$ and let $f,m:
G\rightarrow \mathbb{C}$ be functions such that there exist
functions $M_{1},M_{2}:\rightarrow [0,\infty)$ with
$$\|f(x.y)-f(x)m(y)\|\leq\min\{M_{1}(x),M_{2}(y)\}$$
for all $x, y\in G$. Then either $f$ is bounded or $m(x.y)=m(x)m(y)$
and $f(x)=f(1)g(x)$ for all $x\in G$.
\end{thm}
Also, M. Alimohammady and A. Sadeghi \cite{179} proved a
superstability result for the Cuachy equation. Moreover, they
gived a partial affirmative answer to problem $18$, in the
thirty-first ISFE.

R. Ger pointed out that the superstability phenomenon of the
exponential equation is caused by the fact that the natural group
structure in the range space is disregarded, and he suggested a new
type of stability for the exponential equation (ref. (\cite{ger1}):
\begin{equation}\label{1}
    |\frac{f(x+y)}{f(x)+f(y)}-1|\leq \delta
\end{equation}
If, for each function $f:G\rightarrow E\backslash\{0\}$ satisfying
the inequality (\ref{1}) for some $\delta>0$ and for all $x, y\in
G$, where $E$ is a real Banach space, there exists an exponential
function $M:G\rightarrow E\backslash\{0\}$ such that
\begin{equation}\label{2}
    |\frac{f(x)}{M(x)}-1|\leq \phi(\delta)\ \ \  and\ \ \   |\frac{M(x)}{f(x)}-1|\leq \psi(\delta)
\end{equation}
for all $x\in G$, where $\phi(\delta)$ and $\psi(\delta)$ depend on
$\delta$ only, then the exponential functional equation is said to
be stable in the sense of Ger.

Every complex-valued function of the form $f(x)=a^{x}$ ($x\in
\mathbb{C}$), where $a>0$ is a given number, is a solution of the
functional equation
\begin{equation}\label{2}
    f(xy)=f(x)^{y}
\end{equation}
Hence, the above functional equation may be regarded as a variation
of the exponential functional equation.\cite{jung1}

S.-M. Jung \cite{jun1} proved the stability of the equation
(\ref{2}) in the sense of Ger and he obtained the following result.

\begin{thm}
Let $\delta\in (0,1)$ be a given number. If a function
$f:(0,\infty)\rightarrow (0,\infty)$ satisfies the inequality
\begin{equation}\label{4}
    |\frac{f(xy)}{f(x)^{y}}-1|\leq\delta
\end{equation}
for all $x,y>0$, then there exists a unique constant $a > 0$ such
that $$(1-\delta)^{\alpha(x)}\leq
\frac{a^{x}}{f(x)}\leq(1+\delta)^{\alpha(x)},$$  in which
$\alpha(x)=\sum_{n=1}^{\infty} (\prod_{i=0}^{n-1}x^{2^{i}})^{-1}$
for all $x>1$.
\end{thm}

In this paper, we deal with a type of exponential functional
equation as follows
\begin{equation}\label{3}
   f(xy)=f(x)^{g(y)}
\end{equation}
where $f$ and $g$ are two real valued functions on a commutative
semigroup. We prove that the above functional equation in the sense
of Ger is superstable.

For the reader's convenience and explicit later use, we will recall
a fundamental results in fixed point theory.
\begin{defn}
The pair $(X, d)$ is called a generalized complete metric space if
$X$ is a nonempty set and $d:X^{2}\rightarrow [0,\infty]$ satisfies
the following conditions:
\begin{enumerate}
  \item $d(x,y)\geq0$ and the equality holds if and only if $x=y$;
  \item $d(x,y)=d(y,x)$;
  \item $d(x,z)\leq d(x,y)+d(y,z)$;
  \item every d-Cauchy sequence in X is d-convergent.
\end{enumerate}
for all $x, y\in X$.
\end{defn}
Note that the distance between two points in a generalized metric
space is permitted to be infinity.
\begin{thm}\cite{f4}\label{theo1}
Let $(X, d)$ be a generalized complete metric space and
$J:X\rightarrow X$ be strictly contractive mapping with the
Lipschitz constant $L$. Then for each given element $x\in X$, either
$$d(J^{n}(x),J^{n+1}(x))=\infty$$
for all nonnegative integers $n$ or there exists a positive integer
$n_{0}$ such that
\begin{enumerate}
  \item $d(J^{n}(x),J^{n+1}(x))<\infty$, for all $n\geq n_{0}$;
  \item the sequence $\{J^{n}(x)\}$ converges to a fixed point $y^{\ast}$ of $J$;
  \item $y^{\ast}$ is the unique fixed point of $J$ in the set $Y=\{y\in X\ :\ d(J^{n_{0}}(x),y)<\infty\}$;
  \item $d(y,y^{\ast})\leq\frac{1}{1-L}d(J(y),y)$.
\end{enumerate}
\end{thm}
\section{Main results}
Throughout this Section, assume that $S$ is an arbitrary commutative
semigroup with identity $1$ and $\psi: S^{2}\rightarrow [0,\infty)$
is a function. Let $g:S\rightarrow \mathbb{R}$ be a function, then
we define the set $N_{g}$ with
$$N_{g}=\{a\in S\ :\ |g(a)|>1\}.$$
\begin{thm}\label{theo2}
Suppose that $f:S\rightarrow (0,\infty)$ and $g:S\rightarrow
\mathbb{R}$ are two functions and satisfies the inequality
\begin{equation}\label{4}
    0\leq \frac{f(xy)}{f(y)^{g(x)}}-1\leq \psi(x,y)
\end{equation}
for all $x,y\in S$. If $N_{g}\neq{\O}$ and $\psi(x,ay)\leq\psi(x,y)$
for all $x,y\in S$ and $a\in N_{g}$, then either $g$ is bounded or
\begin{equation}
    f(xy)=f(y)^{g(x)}
\end{equation} for all $x,y\in S$.
\end{thm}
\emph{\textbf{Proof.}} From \ref{4}, we have
\begin{equation}\label{5}
    1\leq \frac{f(xy)}{f(y)^{g(x)}}\leq 1+\psi(x,y)
\end{equation}
for all $x, y\in S$. Then,
\begin{equation}\label{6}
    |\ln \frac{f(xy)}{f(y)^{g(x)}}| \leq
    \ln(1+\psi(x,y))
\end{equation}
or
\begin{equation}\label{7}
     |\ln f(xy) - g(x)\ln(f(y))| \leq \ln(1+\psi(x,y))\leq 1+\psi(x,y)
\end{equation}
for all $x, y\in S$. Set $\tilde{\psi}(x,y):=1+\psi(x,y)$ for all
$x, y\in S$, then its obvious that
\begin{equation}\label{d1}
    \tilde{\psi}(x,ya)\leq\tilde{\psi}(x,y)
\end{equation}
for all $x,y\in S$ and $a\in N_{g}$. Let $a\in N_{g}$ be fixed and
from (\ref{7}), we get
\begin{equation}\label{d2}
     |\ln f(ay) - g(a)\ln(f(y))| \leq \tilde{\psi}(a,y)
\end{equation}
for all $y\in S$. Let us consider the set $A:=\{g:S\rightarrow
(0,\infty)\}$ and introduce the generalized metric on $A$:
$$d(g,h)=\sup_{y\in S}\frac{\|g(y)-h(y)\|}{\tilde{\psi}(a,y)}.$$
It is easy to show that $(A, d)$ is complete metric space. Now we
define the function $J_{a}:A\rightarrow A$ with
$$J_{a}(h(y))=\frac{1}{g(a)}h(ay)$$ for all $h\in A$ and $y\in S$. So
  \begin{eqnarray*}
    d(J_{a}(u),J_{a}(h)) &=& \sup_{y\in S}\frac{\|u(ay)-h(ay)\|}{|g(a)|\tilde{\psi}(a,y)} \\
      &\leq& \sup_{y\in S}\frac{\|u(ay)-h(ay)\|}{|g(a)|\tilde{\psi}(a,ay)}=\frac{1}{|g(a)|}d(u,h)
  \end{eqnarray*}
for all $u, h\in A$, that is $J$ is a strictly contractive
selfmapping of $A$, with the Lipschitz constant
$L=\frac{1}{|g(a)|}$. From (\ref{d2}), we get
$$\|\frac{\ln f(ay)}{g(a)}-\ln f(y)\|\leq\frac{\tilde{\psi}(a,y)}{|g(a)|}$$
for all $y\in S$, which says that $d(J(\ln f),\ln f)\leq L<\infty$.
By Theorem (\ref{theo1}), there exists a mapping $T_{a}:S\rightarrow
(0,\infty)$ such that
  \begin{enumerate}
  \item $T_{a}$ is a fixed point of $J$, i.e.,
  \begin{eqnarray}
    T_{a}(ay) &=& g(a)T_{a}
  \end{eqnarray}
  for all $y\in S$. The mapping $T_{a}$ is a unique fixed point of $J$ in the set $\tilde{A}=\{h\in A\ :\ d(\ln f,h)<\infty\}$.
  \item $d(J^{n}(\ln f),T_{a})\rightarrow 0$ as $n\rightarrow \infty$. This implies that
  $$T_{a}(y)=\lim_{n\rightarrow\infty}\frac{\ln f(a^{n}y)}{g(a)^{n}}$$
  for all $x\in S$.
  \item $d(\ln f,T_{a})\leq\frac{1}{1-L}d(J(\ln f),\ln f)$, which implies,
  $$d(\ln f,T_{a})\leq\frac{1}{|g(a)|-1}.$$
 \end{enumerate}
From (\ref{d2}), its easy to show that following inequality
\begin{eqnarray}\label{d3}
  \|\ln f(a^{n}y)-g(a)^{n}\ln f(y)\| &\leq& \sum_{i=0}^{n-1}\tilde{\psi}(a,ya^{i})|g(a)|^{n-1-i}
\end{eqnarray}
for all $y\in S$ and $n\in \mathbb{N}$. Now since
$\tilde{\psi}(a,ya)\leq\tilde{\psi}(a,y)$ for all $y\in S$, so
$$\tilde{\psi}(a,ya^{m})\leq\tilde{\psi}(a,y)$$ for all $x\in S$ and $m\in
\mathbb{N}$, thus from (\ref{d3}), we obtain
\begin{eqnarray}\label{d4}
  \|\ln f(ya^{n})-g(a)^{n}\ln f(y)\| &\leq& \tilde{\psi}(a,y)\frac{|g(a)|^{n}-1}{|g(a)|-1}
\end{eqnarray}
for all $y\in S$. With this inequality (\ref{d4}), we prove that
$T_{a}=T_{b}$ for each $a, b\in N_{g}$. We have from inequality
(\ref{d4})
\begin{eqnarray}
  \|\ln f(ya^{n})-g(a)^{n}\ln f(y)\| &\leq& \tilde{\psi}(a,y)\frac{|g(a)|^{n}-1}{|g(a)|-1} \label{d5} \\
  \|\ln f(yb^{n})-g(b)^{n}\ln f(y)\| &\leq&\tilde{\psi}(b,y)\frac{|g(b)|^{n}-1}{|g(b)|-1} \label{d6}
\end{eqnarray}
for all $y\in S$. On the replacing $y$ by $yb^{n}$ in (\ref{d5}) and
$y$ by $ya^{n}$ in (\ref{d6})
$$\|\ln f(y(ab)^{n})-g(a)^{n}\ln f(yb^{n})\|\leq \tilde{\psi}(a,y) \frac{|g(a)|^{n}-1}{|g(a)|-1}$$
$$\|\ln f(y(ab)^{n})-g(b)^{n}\ln f(ya^{n})\|\leq\tilde{\psi}(b,y) \frac{|g(b)|^{n}-1}{|g(b)|-1}.$$
Thus,
$$\|g(a)^{n}\ln f(yb^{n})-g(b)^{n}\ln f(ya^{n})\|\leq \tilde{\psi}(a,y) \frac{|g(a)|^{n}-1}{|g(a)|-1}+ \tilde{\psi}(b,y) \frac{|g(b)|^{n}-1}{|g(b)|-1}$$
and dividing by $|g(a)^{n}g(b)^{n}|$
$$\|\frac{\ln f(ya^{n})}{g(a)^{n}}-\frac{\ln f(yb^{n})}{g(b)^{n}}\|\leq$$
$$\frac{\tilde{\psi}(a,y)}{(|g(b)|-1)|g(a)|^{n}}(1-\frac{1}{|g(b)|^{n}})+\frac{\tilde{\psi}(b,y)}{(|g(a)|-1)|g(b)|^{n}}(1-\frac{1}{|g(b)|^{n}})$$
and letting $n$ to infinity, we obtain $T_{a}(y)=T_{b}(y)$ for all
$y\in S$. Therefore, there a unique function $T$ such that $T =
T_{a}$ for every $a\in N_{g}$ and
$$\|\ln f(y)-T(y)\|\leq\frac{\tilde{\psi}(a,y)}{|g(a)|-1}$$ for
all $y\in S$ and $a\in N_{g}$. Since $a\in N_{g}$ is a arbitrary
element, so
$$\|\ln f(y)-T(y)\|\leq\inf_{a\in N_{g}} \frac{\tilde{\psi}(a,y)}{|g(a)|-1} $$
for all $y\in S$. Now if $g$ be a unbounded function, we get $T=\ln
f$.

Let $x, y\in S$ and $a\in N_{g}$ be three arbitrary fixed elements,
from (\ref{d2})
$$\|\ln f(xya^{n})-g(x)\ln f(ya^{n})\|\leq\tilde{\psi}(x,ya^{n})$$ and dividing by $|g(a)|^{n}$,
$$\|\frac{\ln f(xya^{n})}{g(a)^{n}}-g(x)\frac{\ln f(ya^{n})}{g(a)^{n}}\|\leq\frac{\tilde{\psi}(x,ya^{n})}{|g(a)|^{n}}\leq\frac{\tilde{\psi}(x,y)}{|g(a)|^{n}}$$
and letting $n$ to infinity, we get
\begin{equation}
    f(xy)=f(y)^{g(x)}
\end{equation} for all $x,y\in S$. The proof is complete.

With the above Theorem, its easy to obtain the following results.
\begin{cor}
Let $\delta>0$ be a given number. Assume that $f:S\rightarrow
(0,\infty)$ and $g:S\rightarrow \mathbb{R}$ satisfies the inequality
\begin{equation}
    0\leq\frac{f(xy)}{f(y)^{g(x)}}-1\leq\delta
\end{equation}
for all $x,y\in S$, then either $g$ is bounded or
\begin{equation}
    f(xy)=f(y)^{g(x)}
\end{equation}
for all $x,y\in S$.
\end{cor}
\begin{cor}
Let $\delta>0$ be a given number. If a function
$f:(0,\infty)\rightarrow (0,\infty)$ satisfies the inequality
\begin{equation}
    0\leq\frac{f(xy)}{f(y)^{x}}-1\leq\delta
\end{equation}
for all $x,y>0$, then
\begin{equation}
    f(xy)=f(y)^{x}
\end{equation}
for all $x,y>0$.
\end{cor}

\bibliographystyle{amsplain}

\begin{thebibliography}{99}
\bibitem{17} {M. Alimohammady and A. Sadeghi, Stability and common stability for the systems of linear equations and its applications, \em  Math Sci.} 2012, {\bf 6}:43.
\bibitem{179} {M. Alimohammady and A. Sadeghi, Some new results on the superstablity of the Cauchy equation on semigroup, \em Results Math.} (2012) to appear.
\bibitem{Ali} {M. Alimohammady and A. Sadeghi, On the superstability and stability of the Pexiderized exponential equation, \em
CJMS.} {\bf 1}(2)(2012), 61-74.
\bibitem{b2} {J. A. Baker, The stability of the cosine equation, \em Proc. Amer. Math. Soc.} {\bf 80}(1980), 411-416.
\bibitem{b3} {J. A. Baker, J. Lawrence, and F. Zorzitto, The stability of the equation f(x + y) = f(x)f(y), \em Proc. Amer. Math. Soc.} {\bf 74}(1979), 242-246.
\bibitem{asd} {S.Czerwik and M. Przybyla, A general Baker superstability criterion for the D'Alembert functional equation, \em Banach Spaces and their Applications in Analysis,Walter de Gruyter Gmbh @ Co.KG,Berlin} (2007), 303-306.
\bibitem{chung}{J. Chung, On a stability of Pexiderized exponential
equation, \em Bull. Korean Math. Soc.}{\bf 46} (2009), No. 2, pp.
295-301
\bibitem{f4} {J. Diaz and B. Margolis, A fixed point theorem of the alternative for contractions on a generalized
complete metric space, \em Bull. Amer. Math. Soc.} {\bf 74}(1968),
305-309.
\bibitem {ger1}{R. Ger, Superstability is not natural, Rocznik Naukowo-Dydaktyczny WSP w Krakowie, Prace Mat.} 159 (1993), 109-123.
\bibitem{jun1}{S.M. Jung, On the stability of the functional equation $f(xy)=f(x)^{y}$, Mathematica.} (Cluj) 40 (63) (1998), no. 1, 89-94.
\bibitem{jung1} {S.M. Jung, Hyers-Ulam-Rassias stability of functional equations in nonlinear analysis , Springer, Springer Optimization and Its
Appliccations, 2010}
\bibitem{sz1} {L. Sz$\acute{e}$kelyhidi, On a theorem of Baker, Lawrence and Zorzitto, \em  Proc. Amer. Math. Soc.} {\bf 84}(1982), 95-96.
\bibitem{zxc} {J. Tabor and J. Tabor, Homogenity is superstable, \em Publ. Math. Debrecen.} {\bf 45}(1994), 123-130.

\end{thebibliography}

\end{document}